\documentclass[12pt,amsfonts,amssymb,leqno]{article}
\usepackage[english]{babel}
\usepackage[a4paper]{geometry}

\ifx\relax\mathfrak \@xp\@gobbletwo 
  \else \let\mathfrak\relax \fi
\RequirePackage{amsfonts}\relax
\RequirePackage{amsmath}\relax

\RequirePackage{oldlfont}
\def\carre{{$\Box$}}

\makeatletter 
\def\@ypreuve[#1]{\@preuve{ #1}}
\def\@preuve#1{\begin{trivlist}\item[]{\em Proof#1: }} 
\newenvironment{proof}{\@ifnextchar[{\@ypreuve}{\@preuve{}}}{\hfill\carre\end{trivlist}} 

\newtheorem{theorem}{Theorem}
\newtheorem{definition}[theorem]{Definition}
\newtheorem{lemma}[theorem]{Lemma}
\newtheorem{proposition}[theorem]{Proposition}
\newtheorem{corollary}[theorem]{Corollary}
\newtheorem{remark}[theorem]{Remark}

\def\re{{\mathrm{Re}\,}}
\def\im{{\mathrm{Im}\,}}

\def\dom{{\mathrm{dom}}}
\def\ran{{\mathrm{ran}\,}}
\def\ker{{\mathrm{ker}\,}}

\def\dist{{\mathrm{dist}}}

\begin{document}
\begin{center}
{\LARGE\bf Convergence at the origin\\ of integrated semigroups}
\bigskip
\bigskip

by Vincent CACHIA
\footnote{work done at the Department of Theoretical Physics, quai Ernest-Ansermet 24,
CH-1211 GENEVA 4, Switzerland.}

\bigskip
\end{center}

\begin{quote}
{\small
\noindent{\bf Abstract:}
We study a classification of $\kappa$-times integrated semigroups (for $\kappa>0$)
by their (uniform) rate of convergence at the origin: $\|S(t)\|={\cal O}(t^\alpha)$, $t\rightarrow 0$ ($0\leq\alpha\leq\kappa$). By an improved generation theorem we characterize
this behaviour by Hille-Yosida type estimates. Then we consider integrated semigroups with holomorphic
extension and characterize their convergence at the origin, as well as the existence of boundary
values, by estimates of the associated holomorphic semigroup.  Different examples illustrate
these results. The particular case
$\alpha=\kappa$, which corresponds to the notions of Riesz means or tempered
integrated semigroups, is of special interest : as an application, it leads to an integrated
version of Euler's exponential formula.
\medskip

}
\end{quote}
\medskip

\section{Introduction}

\setcounter{equation}{0}
\renewcommand{\theequation}{\arabic{section}.\arabic{equation}}
\setcounter{theorem}{0}
\renewcommand{\thetheorem}{\arabic{section}.\arabic{theorem}}

In the theory of semigroups, convergence at the origin is a key property on which the
standard classification in the well-known treatise of Hille and Phillips \cite{HP} is based.
In the more general theory of integrated semigroups, introduced by Arendt \cite{Arendt},
it seems that no such systematic study has been done. The aim of this paper is to describe
and to characterize the convergence at the origin of integrated semigroups. We cite \cite{ABHN}
as a general reference on the subject. In a first part we consider general integrated semigroups, following \cite{ABHN,Hieber}:
\begin{definition}
Let $A$ be a (multivalued) closed linear operator on a Banach space $X$ with non empty
resolvent set $\rho(A)$, and $\kappa\geq 0$.
The operator $A$ is called the
generator of a $\kappa$-times integrated semigroup if there exist $\omega\geq 0$ and a strongly
continuous function $S:[0,\infty)\rightarrow{\cal L}(X)$ having a Laplace transform for
$\lambda>\omega$, such that $(\omega,\infty)\subset\rho(A)$ and
\begin{equation}
R(\lambda,A) = \lambda^\kappa\int_0^\infty e^{-\lambda t}S(t)dt\quad (\lambda>\omega).
\end{equation}
In this case $S$ is called the $\kappa$-times integrated semigroup generated by $A$.
\end{definition}
The mention of multivalued operators is just a way to include degenerate semigroups
\cite{Yagi}. Hence in general $R(\lambda,A)$ is a pseudo-resolvent, which is not necessarily the
resolvent of a closed univalent operator.
From this definition one deduces that \cite[Proposition 2.4]{Hieber}
\begin{equation}\label{iden1}
S(t)x = \frac{t^\kappa}{\Gamma(\kappa+1)}x + \int_0^t S(s)Axds,\ x\in\dom(A),
\end{equation}
and thus $\lim_{t\rightarrow 0}S(t)x=0$, $x\in\overline{\dom(A)}$, for any integrated semigroup
of order $\kappa>0$. However, at least when $S(t)$ is simply the $\kappa$-th integral of some
$C_0$-semigroup, the expected convergence at the origin is $\|S(t)\|={\cal O}(t^\kappa)$, $t\rightarrow 0$.
Such an estimate has actually been proved for many non-trivial $\kappa$-times
integrated semigroups, that are called tempered (see below \S \ref{fastconvcase} and
\cite{EMK,CarCouOu}).
In fact we show that $\kappa$ is the highest possible power for the convergence at the origin of a non zero $\kappa$-times integrated semigroup ($\kappa>0$). Therefore we aim to study the $\kappa$-times integrated semigroups satisfying $\|S(t)\|={\cal O}(t^\alpha)$, $t\rightarrow 0$, for a given $0\leq\alpha\leq\kappa$.

\section{Hille-Yosida type estimates}

\setcounter{equation}{0}
\renewcommand{\theequation}{\arabic{section}.\arabic{equation}}
\setcounter{theorem}{0}
\renewcommand{\thetheorem}{\arabic{section}.\arabic{theorem}}

Different generation theorems have been proved with Hille-Yosida type estimates, 
cf \cite[Theorem 4.1]{Arendt}, \cite[Theorem 3.4]{Hieber}.
An improvement of these results provides a characterization of the generators of
$\kappa$-times integrated semigroups satisfying $\|S(t)\|={\cal O}(t^\alpha)$, $t\rightarrow 0$, for some $0\leq\alpha\leq\kappa$.

\begin{theorem}\label{carac1}
Let $\{R(\lambda)\}_{\lambda\in\Omega}$ be a pseudo-resolvent on a Banach space $X$ with
$(\omega,\infty)\subset\Omega\subset{\mathbb C}$, for some real $\omega$.
For any $\kappa\geq\alpha\geq 0$
the following are equivalent:

(i) There exist $M>0$ and $a\geq\max\{\omega,0\}$ such that
\begin{equation}\label{HYcondition}
\displaystyle \sup_{n\in{\mathbb N}\cup\{0\},\lambda>a}
\left\|\frac{(\lambda-\omega)^{n+\alpha+1}}{\Gamma(n+\alpha+1)}
\left(\frac{R(\lambda)}{\lambda^\kappa}\right)^{(n)}\right\| \leq M.
\end{equation}

(ii) For each $\delta>0$ there exist a $(\kappa+\delta)$-times integrated semigroup
$\{S_\delta(t)\}_{t\geq 0}$, such that
$R(\lambda)=\lambda^{\kappa+\delta}\int_0^\infty e^{-\lambda t}S_\delta(t)dt$
($\lambda>\omega$), and a constant $M'>0$ such that
\begin{equation}
\|S_\delta(t)-S_\delta(s)\| \leq M'|t-s|^\delta t^{\alpha} e^{\omega t},\ t>s\geq 0.
\end{equation}
\end{theorem}

In fact this theorem appears as a particular case of a result on vector-valued Laplace
transforms (a simpler version of this result can be found in \cite[Corollary 2.5.4]{ABHN}).

\begin{theorem}\label{Laplace}
Let $r\in{\cal C}^\infty((a,\infty),X)$, $a,\alpha\geq 0$ and $\omega\leq a$. Then
the following assertions are equivalent:

(i) There exists $M>0$ such that
\begin{equation}\label{HYlaplace}
\displaystyle \sup_{n\in{\mathbb N}\cup\{0\},\lambda>a}
\left\|\frac{(\lambda-\omega)^{n+\alpha+1}}{\Gamma(n+\alpha+1)}
r^{(n)}(\lambda)\right\| \leq M.
\end{equation}

(ii) For each $\delta>0$ there exists a function
$F_\delta(t):[0,\infty)\rightarrow X$, such that $F_\delta(0)=0$,
$r(\lambda)=\lambda^{\delta}\int_0^\infty e^{-\lambda t}F_\delta(t)dt$
($\lambda>a$), and there exists $M'>0$ such that
\begin{equation}
\|F_\delta(t)-F_\delta(s)\| \leq M'|t-s|^\delta t^{\alpha} e^{\omega t},\ t>s\geq 0.
\end{equation}
\end{theorem}

The proof of this result makes use of fractional integration techniques and includes the
following Lemma.

\begin{lemma}\label{fracinteg}
Let $r\in{\cal C}^\infty((a,\infty),X)$ ($a>0$) and $\delta>0$ such that
$\|\lambda^\delta r(\lambda)\|$ is bounded as $\lambda\rightarrow\infty$.
Then $r_\alpha(\lambda)=\int_\lambda^\infty
\frac{(u-\lambda)^{\alpha-1}}{\Gamma(\alpha)} r(u)du$ is a well defined function 
on $(a,\infty)$ for $0<\alpha<\delta$.
Moreover one has: $(r_\alpha)_\beta=r_{\alpha+\beta}$ with $0<\alpha,\beta<\delta$
and $\alpha+\beta<\delta$,
$(r_\alpha)^{(n)}=(-1)^n r_{\alpha-n}$ with $n\leq\alpha<\delta$
(with the convention $r_0=r$).
If $\sup_{\lambda>a}\lambda^{n+\delta} \|r^{(n)}(\lambda)\|$ is bounded for each $n\in{\mathbb N}$,
then $(r_\alpha)^{(n)}=(r^{(n)})_\alpha$ for all $n\in{\mathbb N}$ and $0<\alpha<\delta$.
If $f$ is a function such that $f_\beta=r_\beta$ for some $0<\beta<[\delta]$, then $f=r$.
\end{lemma}

\begin{proof}
The condition $0<\alpha<\delta$ ensures the convergence of the integral in the sense of Bochner.
Moreover from $\|r(\lambda)\|\leq M/\lambda^\delta$ one deduces the estimate:
\begin{eqnarray*}
\|r_\alpha(\lambda)\| & \leq & \int_\lambda^\infty \frac{(u-\lambda)^{\alpha-1}}
{\Gamma(\alpha)} \frac{M}{u^\delta} du \\
& \leq & \frac{M}{\Gamma(\alpha)}\int_0^1(1-t)^{\alpha-1}\frac{t^{\delta-\alpha-1}}
{\lambda^{\delta-\alpha}}dt\\
& \leq & \frac{M}{\lambda^{\delta-\alpha}}\, \frac{\Gamma(\delta-\alpha)}{\Gamma(\delta)}.
\end{eqnarray*}
If $0<\beta<\delta$ and $\alpha+\beta<\delta$, we have then
\begin{eqnarray*}
(r_\alpha)_\beta(\lambda) & = & \int_\lambda^\infty\frac{(u-\lambda)^{\beta-1}}
{\Gamma(\beta)} \int_u^\infty\frac{(v-u)^{\alpha-1}}{\Gamma(\alpha)} r(v) dudv \\
& = & \int_\lambda^\infty r(v) \int_\lambda^v \frac{(u-\lambda)^{\beta-1}(v-u)^{\alpha-1}}
{\Gamma(\alpha)\Gamma(\beta)} dvdu\\
& = & \int_\lambda^\infty r(v) \frac{(v-\lambda)^{\alpha+\beta-1}}{\Gamma(\alpha+\beta)}dv,
\end{eqnarray*}
which means $(r_\alpha)_\beta=r_{\alpha+\beta}$.

Let $n$ be an integer with $n<\alpha<\delta$, then the integral defining $r_\alpha$
is $n$ times differentiable, which leads to $(r_\alpha)^{(n)}=(-1)^n r_{\alpha-n}$.
If $n=\alpha<\delta$, one obtains $(r_n)^{(n)}=(-1)^n r$.

Let $n\in\mathbb N$ be such that $\sup_{\lambda>a} \lambda^{\delta+n}\|r^{(n)}(\lambda)\|$
is bounded. Then we have $(r^{(n)})_{\alpha+n} = (-1)^n r_\alpha$ by $n$ integrations
by parts. Moreover $[(r^{(n)})_{\alpha+n}]^{(n)}=(-1)^n (r^{(n)})_\alpha$ by the preceding
observation. Finally we obtain $(r_\alpha)^{(n)}=(r^{(n)})_\alpha$.

Suppose that $f_\beta=r_\beta$ for some $\beta$ such that $0<\beta<[\delta]=\sup\{p\in{\mathbb Z}:p\le\delta\}$. Let $\epsilon=[\delta]-\beta$:
then $(f_\beta)_\epsilon=(r_\beta)_\epsilon=r_{[\delta]}$. Moreover since
$[\delta]\in\mathbb N$, one has $r_{[\delta]}^{([\delta])}=(-1)^{[\delta]}r$ and
$f_{[\delta]}^{([\delta])}=(-1)^{[\delta]}f$. Therefore $f=r$.

\end{proof}

\begin{proof}[of Theorem \ref{Laplace}]
If $\alpha=0$, this is the real representation theorem \cite[Theorem 3.2]{Hieber}.
Let now $\alpha>0$, and consider
\begin{equation}
r_\alpha(\lambda)= \int_\lambda^\infty\frac{(u-\lambda)^{\alpha-1}}{\Gamma(\alpha)}r(u)du.
\end{equation}
By Lemma \ref{fracinteg} we have $(r_\alpha)^{(n)}=(r^{(n)})_\alpha$
and we deduce from (\ref{HYlaplace}) that for each $\lambda>a$
\begin{eqnarray*}
\|r_\alpha^{(n)}(\lambda)\| & \leq & \int_\lambda^\infty M\frac{\Gamma(n+\alpha+1)}{\Gamma(\alpha)}
\frac{(u-\lambda)^{\alpha-1}}{(u-\omega)^{n+\alpha+1}}du \\
& \leq & M\frac{\Gamma(n+\alpha+1)}{\Gamma(\alpha)}
         \int_0^1 (1-t)^{\alpha-1}\frac{t^n}{(\lambda-\omega)^{n+1}}dt\\
& \leq & M\frac{n!}{(\lambda-\omega)^{n+1}}.
\end{eqnarray*}
Then for each element $x'$ in the dual space $X'$, by Widder's classical theorem there exists a unique measurable function $f(\cdot,x')$ satisfying $\|e^{-\omega \cdot}f(\cdot,x')\|_\infty
\leq M\|x'\|$ such that
$$
<r_\alpha(\lambda),x'>\, =\int_0^\infty e^{-\lambda t}f(t,x')dt,\ \lambda>a.
$$
In particular one observes that $f(t,x')$ is linear in $x'$.
We set now $\tilde{r}(\lambda)=\int_0^\infty e^{-\lambda t}t^\alpha f(t,x')dt$, $\lambda>a$.
Then we obtain:
\begin{eqnarray*}
\tilde{r}_\alpha(\lambda) &=& \int_\lambda^\infty \frac{(u-\lambda)^{\alpha-1}}{\Gamma(\alpha)}
\int_0^\infty e^{-ut} t^\alpha f(t,x') du dt \\
& = & \int_0^\infty f(t,x')dt\int_{0}^\infty \frac{(vt)^{\alpha-1}}{\Gamma(\alpha)}
e^{-vt-\lambda t}tdv\\
& = & \int_0^\infty e^{-\lambda t} f(t,x')dt=\, <r_\alpha(\lambda),x'>
\end{eqnarray*}
where we have set $v=(u-\lambda)t$. This shows by Lemma \ref{fracinteg} that
$<r(\lambda),x'>\, =\tilde{r}(\lambda)$, $\lambda>a$.
Now for each $\delta>0$ we set $F_\delta(t,x')=\int_0^t \frac{(t-\tau)^{\delta-1}}{\Gamma(\delta)}
\tau^\alpha f(\tau,x')d\tau$ cf \cite[Theorem 6.2.4]{HP}, so that
$$<r(\lambda),x'>\, =\lambda^\delta\int_0^\infty
e^{-\lambda t}F_\delta(t,x')dt.$$
Similarly as in the proof of \cite[Theorem 3.2]{Hieber}, one has
\begin{equation}\label{estHieber}
|F_\delta(t,x')-F_\delta(s,x')|
\leq \frac{2M}{\Gamma(\delta)}|t-s|^\delta t^\alpha e^{\omega t}\|x'\|,\quad t>s\geq 0
\end{equation}
and there exists $F_\delta(t)\in X''$ such that $F_\delta(t,x')=\, <F_\delta(t),x'>$, $x'\in X'$.
Moreover from the uniqueness theorem for the Laplace transform it can be shown that actually
$F_\delta(t)\in X$ for all $t\geq 0$, see \cite{Hieber}.
The norm estimate of $F_\delta(t)$ follows easily from (\ref{estHieber}).

To prove the converse implication, let $\delta=1$ and let $F_1:[0,\infty)\rightarrow X$
be a function which satisfies $(ii)$.
Let $x'\in X'$: since $x'\circ F_1$ is a locally Lipschitz continuous numerical function,
it is differentiable almost everywhere. For $f(t,x')=d/dt<F_1(t),x'>$, we obtain
the estimate $|f(t,x')|\leq M't^\alpha e^{\omega t}\|x'\|$ and
$$
<r(\lambda),x'>\, =\lambda\int_0^\infty e^{\lambda t}<F_1(t),x'>dt
= \int_0^\infty e^{-\lambda t}f(t,x')dt,
$$
for all $\lambda>\max\{\omega,0\}$.
Thus we have
\begin{eqnarray*}
|<r^{(n)}(\lambda),x'>|
& \leq & \int_0^\infty t^n e^{-\lambda t}|f(t,x')|dt \\
& \leq & M'\int_0^\infty t^{n+\alpha}e^{-(\lambda-\omega)t}dt \|x'\| \\
& \leq & M'\frac{\Gamma(n+\alpha+1)}{(\lambda-\omega)^{n+\alpha+1}}\|x'\|,
\end{eqnarray*}
which leads to the assertion $(i)$.
\end{proof}

Now it follows easily that the exponent $\alpha$ in the convergence rate of a non zero $\kappa$-times integrated semigroup cannot be greater than $\kappa$:

\begin{corollary}
Condition $(i)$ of Theorem \ref{carac1} for $\alpha>\kappa$ implies that
$R(\lambda)=0$: this means that a $\kappa$-times integrated semigroup satisfying
$\|S(t)\|={\cal O}(t^\alpha)$ near zero with $\alpha>\kappa$ is necessarily trivial:
$S(t)=0$ for all $t\geq 0$.
\end{corollary}

\begin{proof}
Consider a pseudo-resolvent $\{R(\lambda)\}_{\lambda\in\Omega}$ satisfying (\ref{HYcondition})
 for some $\alpha>\kappa$. It follows (with $n=0$ in (\ref{HYcondition})) that $\sup_{\lambda>a}
\| \lambda^{1+\epsilon}R(\lambda)\| <\infty$ for $\epsilon=\alpha-\kappa>0$.
Then $\sup_{\lambda>a}\| \lambda R(\lambda)\| <\infty$
and by the resolvent equation: $\lim_{\lambda\rightarrow\infty}(\lambda R(\lambda)-I)R(\mu)=0$
for any fixed $\mu$. This shows that $\lim_{\lambda\rightarrow\infty} \lambda R(\lambda)x=x$
for any $x\in\ran R(\mu)$. However $\lim_{\lambda\rightarrow\infty}\|\lambda R(\lambda)\|=0$,
thus the range $\ran R(\mu)$ is reduced to $\{0\}$.
\end{proof}

Arguing as in \cite[Proposition 2.2.2]{ABHN} one can show that a weaker form of the condition
(\ref{HYcondition}) is sufficient:

\begin{corollary}
In order to have the estimate (\ref{HYcondition}), it is sufficient that
\begin{equation}\label{estatn}
\sup_{\lambda>a}\left\|
\frac{(\lambda-\omega)^{n+\alpha+1}}{\Gamma(n+\alpha+1)}
\left(\frac{R(\lambda,A)}{\lambda^\kappa}\right)^{(n)}\right\| \leq M
\end{equation}
for each $n$ in an infinite subset of $\mathbb N$.
Moreover if (\ref{HYcondition}) holds for some $\alpha,\kappa$, then it also holds
for $\alpha+\delta$, $\kappa+\delta$ for any $\delta>0$.
\end{corollary}

If only the issues of growth of integrated semigroups are concerned, then one can state a simpler version of Theorem \ref{carac1}.

\begin{theorem}\label{carac2}
Let $\{S(t)\}_{t\geq 0}$ be a $\kappa$-times integrated semigroup with generator $A$,
$a\geq 0$, $\omega\leq a$ and $0\leq\alpha\leq\kappa$.
Then the following assertions are equivalent:
\begin{eqnarray}
& & \sup_{t>0} t^{-\alpha}e^{-\omega t}\|S(t)\| < \infty, \label{sgt^alpha}\\
& & \sup_{n\in{\mathbb N}\cup\{0\}}\sup_{\lambda>a}\left\|
\frac{(\lambda-\omega)^{n+\alpha+1}}{\Gamma(n+\alpha+1)}
\left(\frac{R(\lambda,A)}{\lambda^\kappa}\right)^{(n)}\right\| < \infty. \label{estHY}
\end{eqnarray}
\end{theorem}

\begin{proof}
Suppose that the $\kappa$-times integrated semigroup $\{S(t)\}_{t\geq 0}$
satisfies (\ref{sgt^alpha}).
Then the associated resolvent is given by
$\frac{R(\lambda,A)}{\lambda^\kappa}=\int_0^\infty e^{-\lambda t}S(t)dt$, and thus
\begin{eqnarray*}
\left\|\left(\frac{R(\lambda,A)}{\lambda^\kappa}\right)^{(n)}\right\|
& \leq & \int_0^\infty Mt^{n+\alpha}e^{-(\re\lambda-\omega)t} dt \\
& \leq & M\frac{\Gamma(n+\alpha+1)}{(\lambda-\omega)^{n+\alpha+1}}
\end{eqnarray*}
for any $\lambda>\max\{\omega,0\}$.

Conversely, suppose that $\{S(t)\}_{t\geq 0}$ is a $\kappa$-times integrated semigroup
and that the resolvent
of the generator satisfies the estimate (\ref{estHY}). Then by the Post-Widder inversion
formula \cite[Theorem 1.7.7]{ABHN} we have for each $x\in X$:
\begin{eqnarray*}
\|S(t)x\| & \leq & \lim_{n\rightarrow\infty}\frac{1}{n!}
\left(\frac{n}{t}\right)^{n+1}\left\|\left(\frac{R(\lambda,A)x}{\lambda^\kappa}
\right)^{(n)}_{\lambda=n/t}\right\| \\
& \leq & \limsup_{n\rightarrow\infty} \frac{\Gamma(n+\alpha+1)}{n!(n/t-\omega)^{n+\alpha+1}}
\left(\frac{n}{t}\right)^{n+1} \|x\| \\
& \leq & M' t^\alpha e^{\omega t}\|x\|.
\end{eqnarray*}

\end{proof}

\section{Euler's exponential formula for $k$-times integrated semigroups}

\setcounter{equation}{0}
\renewcommand{\theequation}{\arabic{section}.\arabic{equation}}
\setcounter{theorem}{0}
\renewcommand{\thetheorem}{\arabic{section}.\arabic{theorem}}

For any once integrated semigroup, an integrated form of Euler's exponential formula
$\lim_{n\rightarrow\infty}(I-tA/n)^{-n}=e^{tA}$
has been established in \cite[Theorem 2.5]{cachia}. We show here a generalization to
$k$-times integrated semigroups ($k\in\mathbb N$), subject to some convergence condition
at the origin.

\begin{theorem}
Let $\{S(t)\}_{t\geq 0}$ be a $k$-times integrated semigroup such that
$\|S(t)\| \leq Mt^\alpha e^{\omega t}$ with
$k-1<\alpha\leq k$ and $k\in\mathbb N$. Let $R(\lambda)=\lambda^k \hat{S}(\lambda)$ be the
associated pseudo-resolvent.
Then the $k$-times integrated Euler's formula holds in the strong operator topology:
define $F(t)=t^{-1} R(t^{-1})$, $0<t<\omega^{-1}$; then for any $t_0>0$,
$F(t/n)^n$ is Bochner integrable on $(0,t_0)$ (for $n$ sufficiently large) and
\begin{equation}
s-\lim_{n\rightarrow\infty} \int_0^{t_0} \frac{(t_0-\tau)^{k-1}}{(k-1)!} F(\tau/n)^n d\tau = S(t_0).
\end{equation}
\end{theorem}

\begin{proof}
Here the integrals make sense as strong Bochner integrals, see e.g. \cite[Ch. 1.1]{ABHN}.
By $R(\lambda)/\lambda^k = \int_0^\infty e^{-\lambda t}S(t)dt$, $\lambda
>\max\{\omega,0\}$, one obtains \cite[Lemma 2.3]{cachia} that
Euler's approximation for $n$ sufficiently large can be written as
\begin{equation}\label{inteuler}
[(n/\tau)R(n/\tau)]^n = \frac{(n/\tau)^{k+1}}{(n-1)!} \int_0^\infty
\left(\frac{nt}{\tau}\right)^{n-k-1} P_{n,k}(nt/\tau) e^{-nt/\tau} S(t)dt,
\end{equation}
where $\tau>0$ and $P_{n,k}$ denotes the polynomial (for $n>k$)
\begin{equation}\label{defPnk}
P_{n,k}(\lambda) = \sum_{\ell=0}^k (-1)^\ell \binom{k}{\ell}
\lambda^{k-\ell} \frac{(n-1)!}{(n-\ell-1)!}.
\end{equation}
One has to verify that Euler's approximation $F(t/n)^n$ is Bochner integrable
near the origin $t=0$. This follows from the convergence condition at the origin
for the $k$-times integrated semigroup. For the term $\ell$ of the sum in (\ref{defPnk})
one has  the estimate (up to some factor depending on $n,\ell$):
\begin{eqnarray*}
& &\hspace{-2cm} \left(\frac{n}{\tau}\right)^{k+1} \int_0^\infty
\left(\frac{nt}{\tau}\right)^{n-\ell-1} e^{-nt/\tau} \|S(t)\|dt\\
& & \leq \frac{M(n/\tau)^{n-\ell+k}}{(n/\tau-\omega)^{n-\ell+\alpha}}
\int_0^\infty u^{n-\ell+\alpha-1} e^{-u}du = {\cal O}\left(\tau^{\alpha-k}\right),
\end{eqnarray*}
as $\tau\rightarrow 0$ (for given $n,\ell$). Thus one has $\|F(\tau/n)^n\| = {\cal O}\left(\tau^{\alpha-k}\right)$ as $\tau\rightarrow 0$ for a given $n$,
and since $\alpha-k>-1$, Euler's approximation is Bochner integrable near the origin.
By using the representation (\ref{inteuler}) and by setting $\lambda=n/\tau$ one has:
\begin{eqnarray}
\lefteqn{ \int_0^{t_0} \frac{(t_0-\tau)^{k-1}}{(k-1)!} F(\tau/n)^n d\tau = }\label{Eulerapp} \\
& & \int_{nt_0^{-1}}^\infty \frac{n d\lambda}{(n-1)!} \int_0^\infty dt S(t)
 \frac{(\lambda t_0-n)^{k-1}}
{(k-1)!} (\lambda t)^{n-k-1} e^{-\lambda t} P_{n,k}(\lambda t).\nonumber
\end{eqnarray}
By Fubini's theorem and the relation
$$
P_{n,k}(\lambda t) \lambda^{n-k-1} e^{-\lambda t} = \frac{d^k}{d\lambda^k}(\lambda^{n-1}
e^{-\lambda t})
$$
\cite[Lemma 2.4]{cachia}, (\ref{Eulerapp}) is equal to
\begin{equation}
\int_0^\infty dt S(t) \frac{nt_0^{k-1}t^{n-k-1}}{(n-1)!} \int_{nt_0^{-1}}^\infty d\lambda
\frac{(\lambda-nt_0^{-1})^{k-1}}{(k-1)!} \frac{d^k}{d\lambda^k}(\lambda^{n-1}e^{-\lambda t}).
\nonumber
\end{equation}
Then by $k-1$ integrations by parts, with all boundary terms vanishing, one finds
\begin{equation}
\int_0^{t_0} \frac{(t_0-\tau)^{k-1}}{(k-1)!} F(\tau/n)^n d\tau =
\frac{n^n}{(n-1)!} \int_0^\infty \frac{dt}{t_0} S(t) \left(\frac{t}{t_0}\right)^{n-k-1} e^{-nt/t_0}.
\nonumber
\end{equation}
Since $S(t)$ is strongly continuous, this integral converges strongly to $S(t_0)$ as
$n\rightarrow\infty$ by an argument similar as in \cite[Theorem 1.3]{cachia} using the Chebycheff inequality.

\end{proof}

\section{Holomorphic extension estimates}

\setcounter{equation}{0}
\renewcommand{\theequation}{\arabic{section}.\arabic{equation}}
\setcounter{theorem}{0}
\renewcommand{\thetheorem}{\arabic{section}.\arabic{theorem}}

In this part we consider integrated semigroups with holomorphic extension into some sector
in the complex plane. We shall study norm estimates in this complex domain and existence of
boundary values, by extending results from \cite{AMH, EMK}.

\subsection{Holomorphic integrated semigroups}

We first state a generation theorem for holomorphic integrated semigroups, and for this
it is convenient to consider also $\kappa$-times integrated semigroup that are not
strongly continuous at the origin:

\begin{definition}
The closed (multivalued) operator $A$ is said to generate a $\kappa$-times
integrated semigroup  $\{S(t)\}_{t>0}$ in the {\bf extended sense} if:
there exists $\omega\geq 0$ such that $(\omega,\infty)\subset\rho(A)$,
$t\mapsto S(t)$ is strongly continuous on $(0,\infty)$,
admits a Laplace transform $\hat{S}(\lambda)$ for $\lambda>\omega$, 
and $\lambda^\kappa\hat{S}(\lambda)=R(\lambda,A)$ for $\lambda>\omega$.
\end{definition}

We mention that if $A$ generates a $\kappa$-times integrated semigroup in the
extended sense, then $A$ generates a $\kappa'$-times integrated semigroup
in the sense of Definition 1.1, for any $\kappa'>\kappa$.

\begin{definition}
A $\kappa$-times integrated semigroup (in the sense of definition 1.1 or 4.1)
is said to be holomorphic of semi-angle $\theta$
if it admits a holomorphic extension into the open sector
$\Sigma_\theta=\{z\in{\mathbb C}: z\neq 0,\ |\arg z|<\theta\}$.
\end{definition}

A direct application of \cite[Theorem 2.6.1]{ABHN} characterises the generators of holomorphic
integrated semigroups in the extended sense:

\begin{theorem}\label{holintegsg}
Let $A$ be a (multivalued) operator in a Banach space $X$, with resolvent set $\rho(A)$.
Then the following assertions are equivalent for $\omega\geq 0$ and $0<\theta\leq\pi/2$:

\begin{tabular}{cp{11cm}}
(i) & The operator $A$ generates a holomorphic $\kappa$-times integrated semigroup of
      semi-angle $\theta$ in the extended sense $S: \Sigma_\theta\rightarrow {\cal L}(X)$
      such that $\sup_{z\in\Sigma_\varphi}\|e^{-\omega z}S(z)\|<\infty$   
      for each $0<\varphi<\theta$.  \\
(ii)& The sector $\omega+\Sigma_{\theta+\pi/2}$ is included in $\rho(A)$ and
      for each $0<\varphi<\theta$
      $$\sup_{\lambda-\omega\in\Sigma_{\varphi+\pi/2}}\|(\lambda-\omega)
      R(\lambda,A)/\lambda^\kappa\|<\infty.$$\\
\end{tabular}\\
Moreover if $(i)$ or $(ii)$ is satisfied, one has the representations:
\begin{equation}\label{integformula}
R(\lambda,A)=\lambda^\kappa\int_0^\infty e^{-\lambda t}S(t)dt\quad
\mbox{ and }\quad S(z)=\frac{1}{2\pi i}\int_C e^{\lambda z}
\frac{R(\lambda,A)}{\lambda^\kappa}d\lambda
\end{equation}
where $C$ is an oriented path in $\omega+\Sigma_{\theta+\pi/2}$,
running from $\infty e^{-i(\varphi+\pi/2)}$ to $\infty e^{i(\varphi+\pi/2)}$ with
$0<\varphi<\theta-|\arg z|$.
\end{theorem}

Let $A$ satisfies one of the assertions $(i)$, $(ii)$ above. For each $\sigma\in\mathbb R$ we consider
the holomorphic function $S_\sigma:\Sigma_\theta\rightarrow{\cal L}(X)$ defined by
\begin{equation}
S_\sigma(z)=\frac{1}{2\pi i}\int_C e^{\lambda z}
\frac{R(\lambda,A)}{\lambda^\sigma}d\lambda.
\end{equation}
Then for $\sigma>\kappa$, $S_\sigma$ is the $\sigma$-times integrated semigroup of generator $A$
(in the sense of Definition 1.1), and $S_{\sigma-k}=S_{\sigma}^{(k)}$ for any $\sigma\in\mathbb R$ and $k\in\mathbb N$. The case $\sigma=0$ is of particular interest, since the holomorphic function $S_0$ satisfies the semigroup equation $S_0(z_1)S_0(z_2)=S_0(z_1+z_2)$ by a standard argument on Cauchy's integrals. This motivates the following definition.

\begin{definition}\label{assholsg}
Let $S_{\kappa}$ be a holomorphic $\kappa$-times integrated semigroup generated by $A$ and of semi-angle $\theta$. We define $S_0(z) = \frac{1}{2\pi i} \int_C e^{\lambda z}R(\lambda,A)d\lambda$ where $C$ is an oriented path in $\omega+\Sigma_{\theta+\pi/2}$,
running from $\infty e^{-i(\varphi+\pi/2)}$ to $\infty e^{i(\varphi+\pi/2)}$ with
$0<\varphi<\theta-|\arg z|$.
The function $S_0:\Sigma_\theta\rightarrow{\cal L}(X)$  is called the holomorphic semigroup associated to the holomorphic $\kappa$-times integrated semigroup $S_{\kappa}$.
\end{definition}
\bigskip

Notice that the function $S_0$ is in general neither bounded, nor Laplace transformable, thus $S_0$ does not fit the standard definition of a holomorphic semigroup. The behaviour of the holomorphic functions $S_\sigma$ (for different values $\sigma$) are related:
we show how to deduce norm estimates and boundary properties of the integrated semigroups $S_\sigma$ (depending on $\sigma$) from a norm estimate of $S_0$ in the complex domain.
For simplicity we consider holomorphic extension in the open right half-plane $\Sigma_{\pi/2}$.

We shall present the results from two points of view. In the first one (section \ref{lths}) the starting point is an operator-valued function on $(0,+\infty)$ satisfying the semigroup equation and admitting a Laplace transform (we then assume that there is a holomorphic extension), whereas in the second one (section \ref{nie}) the starting point is a $\kappa$-times integrated semigroup with holomorphic extension.

\subsection{Laplace transformable holomorphic semigroups}\label{lths}

In this section we consider an operator family $\{T(t)\}_{t>0}\subset{\cal L}(X)$ satisfying the semigroup equation and having a Laplace transform.
This allows to identify the generator $A$ of the semigroup $T$ with the generator of the once integrated semigroup $S_1(t)=\int_0^t T(s)ds$, and we also have $R(\lambda,A)=\hat{T}(\lambda)$. 

Assuming that $T$ admits a holomorphic extension, we show that certain norm estimates for $T(t)$ in the complex domain are equivalent to related convergence estimates at the origin for the integrated semigroups $S_\sigma$. Moreover for $\sigma$ sufficiently large, the holomorphic function $S_\sigma$ admits boundary values on $i\mathbb R$: it extends to a continuous function on the closed right half-plane.

\begin{theorem}\label{Thesthol1}
Let $\{T(t)\}_{t>0}$ be a semigroup on $X$ having a Laplace transform with
abscissa of convergence $\omega_0$.
Then $\hat{T}(\lambda)$ is the resolvent $R(\lambda,A)$ of a
(multivalued) operator $A$ in $X$, with $(\omega_0,\infty)\subset\rho(A)$.
Moreover the following assertions are equivalent for any numbers
$\gamma\geq 0$ and $0\leq\beta<1$:

\begin{tabular}{cp{11cm}}
(i) & The semigroup $T$ has a holomorphic extension into the open right half-plane and
 for each $\alpha>\gamma$ there exist $M>0$ and $\omega\geq 0$ such that
\begin{equation}\label{esthol1}
\|T(z)\| \leq Me^{\omega|z|}\frac{|z|^{\alpha}}{(\re z)^{\alpha+\beta}},\ \re z>0.
\end{equation}\\
(ii) & For each $\alpha>\gamma$, $A$ is the generator of a holomorphic
   $(\alpha+\beta)$-times integrated semigroup $S_{\alpha+\beta}$
   in the open right half-plane, admitting boundary values on $i\mathbb R$,
   and there exist $M,\omega>0$ such that
$\|S_{\alpha+\beta}(z)\| \leq Me^{\omega|z|}|z|^{\alpha}$, $\re z\geq 0$.\\
\end{tabular}\\

\end{theorem}

\begin{proof}
By \cite[Proposition 2.2]{Arendt} the function $\hat{T}$ satisfies the
resolvent equation, thus it is a pseudo-resolvent. It can happen that
$\ker\hat{T}(\lambda)\neq\{0\}$, but in the theory of
multivalued linear operators \cite{Yagi} its inverse is always defined,
and any pseudo-resolvent can be
considered as the resolvent of a multivalued operator. Since $\hat{T}(\lambda)$
is a univalent and bounded operator for $\lambda>\omega_0$, one has
$(\omega_0,\infty)\subset\rho(A)$.

Suppose that $(i)$ is satisfied for some $\gamma\geq 0$ and $0\leq\beta<1$.
Since $\beta<1$, $\|T(z)\|$ is locally integrable at the
origin by estimate (\ref{esthol1}).
Then we follow an argument due to El-Mennaoui, cf \cite[Theorem 5.1]{CarCouOu}.
For each $\alpha'>0$ the $\alpha'$-times integrated semigroup generated by $A$ is given by
\begin{equation}
	 S_{\alpha'}(z)=\frac{1}{\Gamma(\alpha')}\int_0^z (z-s)^{\alpha'-1}T(s)ds,\ \re z>0, \label{Salpha}
\end{equation}
where the integral is absolutely convergent
and does not depend on the path from $0$ to $z$ in the open right half-plane.
Without loss of generality, let $\arg(z)>0$.
We then consider the path consisting of two parts: the straight line $[0,|z|]$ and the arc
$\{|z|e^{it}: 0\leq t\leq\arg z\}$. We call $S_{\alpha'}^1(z)$ and $S_{\alpha'}^2(z)$ the integrals
corresponding to each part, and estimate them separately with the help of (\ref{esthol1}),
for some $\alpha>\gamma$. For the first part we have
$$
\|S_{\alpha'}^1(z)\| \leq
\frac{Me^{\omega|z|}}{\Gamma(\alpha')}\int_0^{|z|} \frac{|z-s|^{\alpha'-1}}{s^{\beta}}ds
\leq M e^{\omega|z|}|z|^{\alpha'-\beta}\frac{C_{\alpha',\beta}}{\Gamma(\alpha')},
$$
where $C_{\alpha',\beta}=\sup_{|\theta|\leq\pi/2}\int_0^1\frac{|e^{i\theta}-u|^{\alpha'-1}}
{u^{\beta}}du<\infty$.
For the second part we have:
$$
\|S_{\alpha'}^2(z)\| \leq \frac{Me^{\omega|z|}}{\Gamma(\alpha')}\int_0^{\arg z}
|z|^{\alpha'-1}|e^{i\arg z}-e^{it}|^{\alpha'-1}\frac{|z|^\alpha}{(|z|\cos t)^{\alpha+\beta}} |z|dt.
$$
By observing that $|e^{i\arg z}-e^{it}|=2\sin(\frac{\arg z-t}{2})$ and that $\cos t\ge \sin(\frac{\arg z-t}{2})$,
and by setting $\theta=(\arg z-t)/2$ we obtain the estimate:
\begin{equation}
\|S_{\alpha'}^2(z)\| \leq\frac{Me^{\omega|z|}}{\Gamma(\alpha')}
|z|^{\alpha'-\beta}\int_0^{\frac{\arg z}{2}}
2^{\alpha'}(\sin\theta)^{\alpha'-\alpha-\beta-1}d\theta,\label{Salpha2}
\end{equation}
where the last integral is finite provided $\alpha'>\alpha+\beta$. Therefore for each
$\alpha'>\gamma+\beta$ one can choose $\alpha>\gamma$ such that $\alpha'>\alpha+\beta$ and
then $\|S_{\alpha'}(z)\| \leq M' e^{\omega|z|}|z|^{\alpha'-\beta}$,
$\re z>0$.
Moreover the integral in (\ref{Salpha2}) also converges for $\arg(z)=\pi/2$, hence the integral (\ref{Salpha}) with $z\in i\mathbb R$ defines bounded boundary values for $S_{\alpha'}$ with the same estimate. Thus assertion $(ii)$ is proved.
\bigskip

Conversely suppose that $(ii)$ is satisfied, and let $\alpha>\gamma$. For any integer $k>\alpha+\beta$
we obtain by integration of $S_{\alpha+\beta}$, cf (\ref{estHieber})~:
\begin{equation}\label{estS_k}
\|S_k(z+h)-S_k(z)\| \leq M(|z|+|h|)^\alpha |h|^{k-\alpha-\beta}e^{\omega(|z|+|h|)}.
\end{equation}
Then by Cauchy's formula
$$
T(z)=\frac{k!}{2\pi i}\int_{C_z}\frac{S_k(\zeta)}{(\zeta-z)^{k+1}}d\zeta
=\frac{k!}{2\pi i}\int_{C_z}\frac{S_k(\zeta)-S_k(z)}{(\zeta-z)^{k+1}}d\zeta
$$
where $C_z$ denotes the path defined by the circle with centre $z$ and radius $r=\frac{\re z}{2}$.
Then by (\ref{estS_k}) we have
\begin{eqnarray*}
\|T(z)\| & \leq & \frac{k!}{2\pi r^{k}}\left\|\int_0^{2\pi}(S_k(z+re^{it})-S_k(z))e^{-ikt}dt\right\|\\
& \leq &\frac{k!}{r^{\alpha+\beta}}M \left|\frac{3z}{2}\right|^\alpha e^{3\omega|z|/2}
\leq k!2^\beta 3^\alpha M\frac{|z|^\alpha}{(\re z)^{\alpha+\beta}} e^{3\omega|z|/2},
\end{eqnarray*}
which leads to $(i)$.

\end{proof}

\begin{remark}
In order that a strongly continuous (for $t>0$) semigroup $\{T(t)\}_{t>0}$
admits a Laplace transform, it is sufficient that
$\int_0^1 \|T(t)\|dt<\infty$. This follows from the observation:
$\int_n^{n+1} \|T(t)\|dt \leq \|T(1)\|^n \int_0^1 \|T(t)\|dt$.
In fact for each $t_0>0$, $\|T(t)\|$ is exponentially bounded for $t\geq t_0$.
\end{remark}

\begin{remark}\label{extesthol}
Let $\{T(z)\}_{\re z>0}$ be a holomorphic semigroup and $\alpha,\beta\geq 0$.
In order that there exist $\omega,M>0$ such that the estimate (\ref{esthol1}) holds
in the open right half-plane,
it is sufficient that such an estimate holds in the vertical strip $0<\re z<1$.
\end{remark}

\begin{proof} 
Let us suppose that (\ref{esthol1}) holds for some $\omega_0,M_0>0$ and
$0<\re z<1$. Then we have for some $M\geq M_0$ and $\omega\geq \omega_0$:
$$
\|T(z)\| \leq M e^{\omega|z|},\ 1/2<\re z<1.
$$
For any number $z=a+ib$ such that $a>1$ we set $k=[a]+1\in{\mathbb N}$ so that
$1/2\leq a/k<1$. Therefore $1/2\leq\re z/k<1$ and we have:
$$
\|T(z)\| \leq \|T(z/k)\|^k \leq M^k e^{\omega|z|}\leq M e^{(\omega+\ln M)|z|},\ \re z>1.
$$
This leads to an estimate of type (\ref{esthol1}) in the open right half-plane, with possibly new constants $M,\omega$.

\end{proof}

\subsection{Non integrable estimates}\label{nie}

The question arises: what could be deduced from an estimate of type (\ref{esthol1}) with $\beta\geq 1$ (which does not ensure the integrability near the origin) ? 
In this case we are not, in general, able to express the resolvent of the generator by Laplace transform of the semigroup. To circumvent this difficulty we assume that the semigroup $T(z)$ actually is the derivative of some $\kappa$-times integrated semigroup.
The following result is an improvement of \cite[Theorem 3.9.13]{ABHN}, on the existence
of boundary values of holomorphic integrated semigroups.

\begin{theorem}\label{Thesthol2}
Let $A$ generate a holomorphic $\kappa$-times integrated semigroup $S_\kappa$
in the open right half-plane
for some $\kappa>0$, and let $T=S_0$ denote the holomorphic semigroup associated to $S_\kappa$ by Definition \ref{assholsg}.
Let us consider the following assertions for $\alpha,\beta,\gamma,\delta\geq 0$
\begin{eqnarray*}
& (E^1_{\alpha,\beta}) & \exists\omega_1,M_1>0 \mbox{ such that } \|T(z)\|\leq M_1
 e^{\omega_1|z|}\frac{|z|^\alpha}{(\re z)^{\alpha+\beta}},\ \re z>0.\\
& (E^2_{\gamma,\delta}) & \exists\omega_2,M_2>0 \mbox{ such that }
 \|S_{\gamma+\delta}(z)\| \leq M_2 e^{\omega_2|z|}|z|^{\gamma},\ \re z\geq 0.
\end{eqnarray*}
Then one has $(E^1_{\alpha,\beta})\Rightarrow (E^2_{\gamma,\delta})$ for each $\gamma>\alpha$
and $\delta>\beta$, and $(E^2_{\gamma,\delta}) \Rightarrow (E^1_{\gamma,\delta})$.
\end{theorem}

\begin{proof}
If $\beta<1$, then the result follows from Theorem \ref{Thesthol1}.
Now we suppose that $\beta\geq 1$.
By \cite[Proposition 3.2.6]{ABHN}, we may rescale the problem without loss
of generality and consider $A-a$ instead of $A$ ($a>0$). Then the $k$-times integrated semigroup
generated by $A-a$ is related to that generated by $A$ as follows:
$$
S_a(t) = e^{-at}S(t)+\sum_{j=1}^k \binom{k}{j}a^j\int_0^t\frac{(t-s)^{j-1}}{(j-1)!}
e^{-as}S(s)ds.
$$
Thus the estimate $\|S(t)\| \leq Mt^\alpha e^{\omega t}$ implies
$\|S_a(t)\| \leq M_a t^\alpha e^{\omega_a t}$ for some $M_a\geq M$ and $\omega_a\leq\omega$.
This shows that the assertion $(E^2_{\gamma,\delta})$ for $S_a$ is equivalent to that for $S$.
For the associated holomorphic semigroup, $T(z)$ is replaced by $T_a(z)=e^{-az}T(z)$ by the
integral formula (\ref{integformula}). Thus for $a>\omega$ the estimate
$(E^1_{\alpha,\beta})$ gives:
\begin{equation}\label{esthol1'}
\|T_a(z)\| \leq Me^{\omega_i|\im z|}\frac{|z|^\alpha e^{-\omega_r\re z}}{(\re z)^{\alpha+\beta}}
\end{equation}
for some $M,\omega_i,\omega_r>0$, which is again of the same type. From now on we consider
the rescaled problem and omit the mention of $a$.

For any integer $k$ such that $A$ generates a $k$-times integrated semigroup $S_k$ one has by
iteration of (\ref{iden1}):
\begin{equation}\label{dkSk}
\frac{d^k}{dt^k} S_k(t)x=S_k(t)A^k x + \frac{t^{k-1}}{(k-1)!}A^{k-1}x + \dots + tAx+x,\ x\in\dom(A^k).
\end{equation}
Since by hypothesis $A$ generates a $\kappa$-times integrated semigroup,
the identity (\ref{dkSk}) holds for any integer $k\geq\kappa$. Hence $T(t)x=\frac{d^k}{dt^k}S_k(t)x$
is Laplace transformable for $x\in\dom(A^k)$ and $k\geq\kappa$,
and we have for each $x\in\dom(A^k)$
$$
R(\lambda,A)x = \int_0^\infty e^{-\lambda u}T(u)x du,\ \lambda>-\omega_r,
$$
which leads to
$$
R(\lambda,A)^{n} x = \frac{(-1)^{n-1}}{(n-1)!} R(\lambda,A)^{(n-1)}x =
\int_0^\infty e^{-\lambda u}\frac{u^{n-1}}{(n-1)!}T(u)xdu,\ n=1,2,\dots
$$
Let us consider an arbitrary $z=t+is$ with $\re z=t>0$: since $T(z)$ is a holomorphic semigroup,
$\ran T(z)\subset\dom(A^k)$ for each $k\in\mathbb N$, in particular for $k\geq\kappa$.
Then we have for all $x\in X$
$$
A^{-n}T(z)x = \int_0^\infty \frac{u^{n-1}}{(n-1)!} T(u+z)x du,\ n=1,2,\dots
$$
Then by the estimate (\ref{esthol1'}) one obtains 
$$
\|A^{-n}T(z)x\| \leq \frac{Me^{\omega_i|s|}}{(n-1)!} \int_0^\infty u^{n-1}
\frac{|u+z|^\alpha}{(u+t)^\alpha}  \frac{e^{-\omega_r (u+t)}}{(u+t)^{\beta}}\|x\| du,\ n=1,2,\dots
$$
We observe that $1<\frac{|u+z|}{\re(u+z)}\leq\frac{|z|}{\re z}$ for all $u\geq 0$.
Let us now choose $n=[\beta]=\sup\{p\in{\mathbb Z}:p\le\beta\}\geq 1$ and a number $\epsilon\in(0,1)$ and let us set $\beta'=\beta-[\beta]+\epsilon$. By observing that
$(u+t)^{-\beta} = (u+t)^{-[\beta]+\epsilon} (u+t)^{-\beta'} \leq u^{-[\beta]+\epsilon} t^{-\beta'}$ for all $u,t>0$,
one finds for all $x\in X$
\begin{eqnarray*}
\|A^{-[\beta]}T(z)x\| & \leq & \frac{M|z|^\alpha e^{\omega_i|s|}}{(\re z)^\alpha([\beta]-1)!}
 \int_0^\infty
u^{\epsilon-1} \frac{e^{-\omega_r (u+t)}}{t^{\beta'}}\|x\| du\\
& \leq & M'\frac{|z|^\alpha e^{\omega_i|s|}e^{-\omega_r t}}{t^{\alpha+\beta'}}\|x\|.
\end{eqnarray*}
By (\ref{dkSk}) one deduces that for each
$\beta'>\beta-[\beta]$, there exist $\tilde{M},\omega>0$ such that
$$
\|S_{[\beta]}(z)\| \leq \frac{\tilde{M}|z|^\alpha e^{\omega|z|}}{(\re z)^{\alpha+\beta'}}.
$$
By choosing $\beta'<1$ it follows that the integral
$$
\int_0^z \frac{(z-\zeta)^{\alpha'-1}}{\Gamma(\alpha')} S_{[\beta]}(\zeta)d\zeta
$$
converges for $\re z>0$ and $\alpha'>0$; moreover it clearly coincides with the holomorphic
$([\beta]+\alpha')$-times integrated semigroup generated by $A$.
By using an estimation method similar to that used for the proof of Theorem \ref{Thesthol1},
one obtains that for each $\alpha'>\alpha+\beta'$, $S_{[\beta]+\alpha'}$ admits boundary
values on $i\mathbb R$ and there exist $M_0,\omega>0$ such that
$$
\|S_{[\beta]+\alpha'}\| \leq M_0 |z|^{\alpha'-\beta'} e^{\omega|z|},\ \re z\geq 0,
$$
which is $(E^2_{\gamma,\delta})$ by setting $\gamma=\alpha'-\beta'>\alpha$ and
$\delta=[\beta]+\beta'>\beta$.
\bigskip

The converse implication $(E^2_{\gamma,\delta})\Rightarrow(E^1_{\gamma,\delta})$
is proved exactly as in Theorem~\ref{Thesthol1}.
\end{proof}

We now complete the study by establishing estimates for all functions $S_\sigma$, similar to 
 $(E^2_{\gamma,\delta})$.
We mention that a holomorphic $(\gamma+\delta)$-times integrated semigroup satisfying
the estimate $(E^2_{\gamma,\delta})$ is always obtained by integration
of some holomorphic $\delta$-times integrated semigroup in the extended sense.
This property is not satisfied by integrated semigroups without holomorphic extension
(see Proposition \ref{hieber}).

\begin{corollary}
Let $A$ generate a holomorphic $\kappa$-times integrated semigroup in the open right
half-plane and suppose that estimate $(E^1_{\alpha,\beta})$ holds.
Then for each $\gamma>\alpha$, $\delta>\beta$ and $\sigma<\gamma+\delta$ we have
$$
(E^3_{\sigma,\gamma,\delta})\quad\exists\omega_3,M_3>0 \mbox{ such that }
 \|S_{\sigma}(z)\| \leq M_3
e^{\omega_3 |z|} \frac{|z|^{\gamma}}{(\re z)^{\gamma+\delta-\sigma}},\ \re z>0.
$$

For each $\sigma>\beta$, $S_\sigma$ is a holomorphic $\sigma$-times integrated semigroup of generator
$A$, with exponential bounds in each sector $\Sigma_\theta$,
$\theta<\pi/2$.

For each $\sigma>\alpha+\beta$, $iA$ is the generator of a $\sigma$-times integrated group:
$e^{-i\sigma\pi/2}S_\sigma(it)$ for $t>0$ and $e^{i\sigma\pi/2}S_\sigma(it)$ for $t<0$,
with bound $Me^{\omega|t|}|t|^{\sigma-\delta}$, $t\in\mathbb R$, for each $\delta>\beta$.
\end{corollary}

\begin{proof}
By Theorem \ref{Thesthol2}, the estimate $(E^1_{\alpha,\beta})$ implies $(E^2_{\gamma,\delta})$
for any $\gamma>\alpha$ and $\delta>\beta$.
Similarly to the implication $(ii)\Rightarrow(i)$ in Theorem \ref{Thesthol1}, we obtain that
$(E^2_{\gamma,\delta})\Rightarrow(E^3_{\sigma,\gamma,\delta})$.
Then it follows that $S_\sigma$ is exponentially bounded in each sector $\Sigma_\theta$,
$\theta<\pi/2$, for $\sigma>\beta$.

If $\sigma>\alpha+\beta$, then one can choose $\gamma>\alpha$ and $\delta>\beta$ such that
$\sigma=\gamma+\delta$ and $S_\sigma$ satisfies $(E^2_{\gamma,\delta})$. Now we show that
\begin{equation}\label{integtheta}
\int_0^\infty e^{-\lambda t}S_\sigma(t)dt =
\int_{(0,\infty e^{i\theta})} e^{-\lambda z}S_\sigma(z)dz
\end{equation}
for $0\leq\theta\leq\pi/2$ and $\lambda$ in some domain to be specified.
Let us close the integration path consisting of the segments $[0,R]$, $[0,Re^{i\theta}]$
by the arc $\Gamma_R=\{|z|e^{it} : 0\le t\le \theta\}$. We have to estimate the integral
on $\Gamma_R$:
$$
\left\|\int_{\Gamma_R} e^{-\lambda z}S_\sigma(z)dz\right\| \leq \int_0^\theta e^{-\re(\lambda
 Re^{i\varphi})}\|S_\sigma(Re^{i\varphi})\|Rd\varphi.
$$
Let $\eta\in(0,\pi/2)$: 
then one has $\re(\lambda e^{i\varphi})\geq |\lambda|\cos\eta$
for any $0\leq\varphi\leq\theta$ and for any $\lambda$ such that
$-\pi/2+\eta\leq\arg\lambda\leq\pi/2-\theta-\eta$. Therefore for such $\lambda$
the integral on $\Gamma_R$ is bounded by
$M_2 \theta R^{\gamma+1}e^{R(\omega_2 -|\lambda|\cos\eta)}$, and thus tends to $0$ as
$R\rightarrow\infty$ provided $|\lambda|$ is sufficiently large.
This leads to the equality (\ref{integtheta}) for these specified values of $\lambda$.

Since $\int_0^\infty e^{-\lambda t}S_\sigma(t)dt$ coincides with the holomorphic function
$R(\lambda,A)/\lambda^\sigma$ (for $\re\lambda>\omega_0$, the abscissa of convergence),
one deduces that
$$
\lambda^\sigma\int_0^\infty e^{-\lambda e^{i\theta}t}S_\sigma(e^{i\theta}t)e^{i\theta}dt
= R(\lambda,A) = e^{i\theta} R(\lambda e^{i\theta},Ae^{i\theta}).
$$
for any $\lambda$ such that the integral is defined. Hence
$$
R(\mu,e^{i\theta}A) = e^{-i\sigma\theta}\mu^\sigma\int_0^\infty e^{-\mu t}
S_\sigma(e^{i\theta}t)dt,
$$
which means that $\{e^{-i\sigma\theta}S_\sigma(e^{i\theta}t)\}_{t\geq 0}$ is a $\sigma$-times
integrated semigroup with the generator $e^{i\theta}A$. Since the same is true for $\theta=\pi/2$
and $\theta=-\pi/2$ one obtains that $iA$ generates a $\sigma$-times integrated group.

\end{proof}

\begin{remark}
In fact, Theorems \ref{Thesthol1} or \ref{Thesthol2}
could be stated in the more general
case where the analyticity sector $\Sigma_\theta$ has a semi-angle $\theta\leq\pi/2$:
in the estimates one has just to replace $\re z$ by $\dist(z,\partial\Sigma_\theta)$,
where the boundary $\partial\Sigma_\theta$ is $e^{\pm i\theta}{\mathbb R}_+$.
\end{remark}

\section{Examples}

\setcounter{equation}{0}
\renewcommand{\theequation}{\arabic{section}.\arabic{equation}}
\setcounter{theorem}{0}
\renewcommand{\thetheorem}{\arabic{section}.\arabic{theorem}}

\subsection{The optimal convergence case $\alpha=\kappa$}\label{fastconvcase}

Many examples of integrated semigroups have been shown to exhibit such behaviour.
These integrated semigroups are also called tempered, and they admit bounded
Riesz means. The typical example is the Gaussian semigroup $T_p(z)=e^{z\Delta_p}$
in $L^p({\mathbb R}^n)$, which satisfies the estimates \cite[(3.61) and (3.62)]{ABHN}
$$ 2^{-n/2p}\left(\frac{|z|}{\re z}\right)^{n|1/p-1/2|} \leq
\|T_p(z)\|_{{\cal L}(L^p({\mathbb R}^n))}
\leq \left(\frac{|z|}{\re z}\right)^{n|1/p-1/2|}$$
for $\re z>0$ and $1\leq p\leq\infty$.
For the Schr\"odinger group, similar results were obtained in \cite{CarCouOu}.
We mention here a particularly interesting result, which is taken from
{\cite[Theorem 4.2, 4.3]{Hieber2}}: it shows that the convergence at the origin is optimal
for a large class of differential operators. Moreover the associated $\kappa$-times integrated semigroup admits boundary values only for sufficiently large $\kappa$.

\begin{proposition}\label{hieber}
let $A$ be a differential operator on $L^p({\mathbb R}^n)$ ($1\leq p\leq\infty$) with maximal
domain such that its symbol is of the form $ia(\xi)$, where $a(\xi)$ is a real, homogeneous,
elliptic polynomial on ${\mathbb R}^n$.
Then $A$ generates an $\alpha$-times integrated semigroup $S(t)$ on $L^p({\mathbb R}^n)$
satisfying $\|S(t)\| \leq Mt^\alpha$ for some $M>0$ and all $t>0$, whenever $\alpha>n|1/2-1/p|$.
Moreover if the symbol of $A$ is of the form $\pm i|\xi|^m$ for some $m>0$ with $m\neq 1$,
then $A$ generates an $\alpha$-times integrated semigroup: if and only if
$\alpha>n/2$ for $p=1$ or $p=\infty$, if and only if $\alpha\geq n|1/2-1/p|$ for
$1<p<\infty$.
\end{proposition}

In particular, for any $\kappa>0$
there is a generator of a $\kappa$-times integrated semigroup with optimal convergence
rate, which does not generate a $\sigma$-times integrated semigroup for $\sigma<\kappa$.

Another situation where the optimal convergence takes place is the fractional powers
of an operator. Let $B$ be an operator in $X$ such that $(-\infty,0]\subset\rho(B)$
and $\sup_{\lambda\leq 0}(1-\lambda)\|R(\lambda,B)\|<\infty$. Then (following
\cite[Th. 3.8.1]{ABHN}) the fractional power
$$
B^{-z} = \frac{1}{2\pi i}\int_\Gamma \lambda^{-z} R(\lambda,B)d\lambda,\ \re z>0,
$$
(where $\Gamma$ is a smooth path in $\rho(B)$ going from $\infty e^{-i\delta}$ to $\infty e^{i\delta}$ for some $\delta>0$) defines a holomorphic semigroup on $X$, and satisfies the estimate:
$$
\|B^{-z}\| \leq M\frac{|\sin\pi z|}{\sin(\pi\re z)},\ 0<\re z<1.
$$
We thus have an estimate of type $e^{\omega|z|}|z|/\re z$ (cf Remark \ref{extesthol}).
By Theorem \ref{Thesthol1}, for each $\kappa>1$, the associated $\kappa$-times integrated semigroup admits boundary values, which satisfy an estimate near the origin with the optimal exponent $\kappa$.

The property $\alpha=\kappa$ is a useful hypothesis for different other
results on integrated semigroups, e.g. \cite[Theorem 4.1]{Hieber}.

\subsection{The intermediate case $0\leq\alpha\leq\kappa$}

It seems that no example outside the optimal convergence case
has been explicitly described.
The following one is inspired by \cite[\S I.8.1]{Krein}, and illustrates all different
cases by changing a parameter $\beta$.
Let us consider the Banach space $X=L^p({\mathbb R})\times L^p({\mathbb R})$
for some $p\in[1,\infty]$, with norm $\|(u,v)\| = {\|u\|}_p + {\|v\|}_p$,
and $\beta\geq 0$. We define a multiplication
operator $A$ on $X$ by the matrix:
\begin{equation}
a(x)=\left(\begin{array}{cc}
-(1+x^{2}) & |x|^\beta \\
0 & -(1+x^{2}) \\
\end{array}\right) = -(1+x^{2})I + |x|^\beta N,
\end{equation}
where $x\in\mathbb R$, $I$ denotes the identity $2\times 2$ matrix, and
$N=\left(\begin{array}{cc} 0 & 1 \\ 0 & 0 \\ \end{array} \right)$.
It follows that $a(x)^n = (-1-x^{2})^n I + n(-1-x^{2})^{n-1} |x|^\beta N$ and then
$e^{ta(x)} = e^{-t(1+x^{2})} (I+t|x|^\beta N)$ for each $x\in\mathbb R$ and $t\in\mathbb C$.

\begin{lemma}
The multiplication operator $T(t)$ associated to $e^{ta(x)}$ is bounded on $X$ for $\re t>0$
and $\beta\geq 0$.
The operator-valued function $t\mapsto T(t)$ is holomorphic in the open right
half-plane and satisfies the semigroup equation $T(t+s)=T(t)T(s)$.
Moreover one has the estimate:
$$\max\left\{1,C_\beta\frac{|t|e^{-\re t}}{(\re t)^{\beta/2}}\right\}\leq \|T(t)\| \leq
1+C_\beta\frac{|t|e^{-\re t}}{(\re t)^{\beta/2}}$$
where $C_\beta=(\beta/2)^{\beta/2}e^{-\beta/2}$.

If $0\leq\beta<2$, $\{T(t)\}_{\re t>0}$ is a holomorphic $C_0$-semigroup,
if $2\leq\beta<4$, it is of class $(1,A)$,
if $\beta=4$ it is still Abel summable,
and finally if $\beta>4$ the resolvent set of $A$ is empty.
\end{lemma}

\begin{proof}
The argument is based on the following estimate for $\re t>0$:
\begin{equation}\label{estbetak}
\sup_{x\in\mathbb R} |x|^\beta |te^{-t(1+x^{2})}| = \sup_{x\in\mathbb R}
\frac{|t|e^{-\re t}}{(\re t)^{\beta/2}} (x^{2}\re t)^{\beta/2} e^{-x^{2}\re t}
= C_\beta \frac{|t|e^{-\re t}}{(\re t)^{\beta/2}},
\end{equation}
where $C_\beta=\sup_{y>0}y^{\beta/2} e^{-y}=(\beta/2)^{\beta/2}e^{-\beta/2}$.
This gives lower and upper bounds for the
norm of the multiplication operator $e^{ta(x)}$. If $0\leq\beta<2$,
the convergence at $t\rightarrow 0$ is easily verified and we have a $C_0$-semigroup.

From the expression of $a(x)$ we obtain that
$(\lambda-A)^{-1}$ is the multiplication operator by the matrix:
$(\lambda-a(x))^{-1}=(\lambda+1+x^2)^{-1}I + |x|^\beta(\lambda+1+x^2)^{-2}N$.
This shows that $\rho(A)={\mathbb C}\setminus(-\infty,-1]$ for $0\leq\beta\leq 4$
and $\rho(A)=\emptyset$ for $\beta>4$. Furthermore one has for $\beta\leq 4$:
$\lim_{\lambda\rightarrow\infty}\lambda R(\lambda,A)f=f$ for each $f\in X$,
which means that $T(t)$ is Abel summable \cite[\S10.6]{HP}. When $\beta<4$ it is
also integrable near the origin, thus it is in the class $(1,A)$.
\end{proof}

The semigroup $\{T(t)\}_{\re t>0}$ satisfies the estimate $(E^1_{1,\beta/2-1})$.
If $t\notin\mathbb R$ and $\beta>0$, the operator-norm of $T(t)$ becomes infinite as $\re t\rightarrow 0$,
which shows that $\{T(t)\}_{\re t>0}$ has no boundary values on $i\mathbb R$.
However Theorem~\ref{Thesthol1} applies for $\beta<4$ and shows that
the associated $\kappa$-times integrated semigroup admits boundary values for sufficiently large $\kappa$.

\begin{lemma}
Let $0\leq\beta\leq 4$.
For each $\kappa>\beta/2$, the operator $A$ generates
a holomorphic $\kappa$-times integrated semigroup with boundary values on $i\mathbb R$.
Moreover one has the estimates (for some $M>0$):
$\|S_\kappa(z)\| \leq M|z|^{\kappa}$ if $0\leq\beta\leq 2$, and
$\|S_\kappa(z)\| \leq M|z|^{1+\kappa-\beta/2}$ if $2\leq\beta\leq 4$.
\end{lemma}

\begin{proof}
For $\beta<4$ the result follows directly from Theorem \ref{Thesthol1}.
The case $\beta=4$ requires more attention.
In fact, an explicit calculation gives the once integrated semigroup of generator $A$:
for $t>0$, $S_1(t)$ is the multiplication operator associated to the matrix:
\begin{equation}\nonumber
s_1(x,t)=\int_0^t e^{\tau a(x)}d\tau = \frac{1-e^{-t(1+x^{2})}}{1+x^{2}} I
+ |x|^\beta N \left[\frac{-t e^{-t(1+x^{2})}}{1+x^{2}} +
\frac{1-e^{-t(1+x^{2})}}{(1+x^{2})^2}\right].
\end{equation}
Hence by (\ref{estbetak}), $S_1(t)$ is a bounded operator for $t>0$ and $0\leq\beta\leq 4$.
Then $S_1(t)$ has a bounded holomorphic extension to the open right half-plane, which admits
boundary values for $t\in i\mathbb R$, whenever $0\leq\beta\leq 2$.
Moreover by observing that for $\re t\geq 0$,
$(1+x^{2})^{-1}|1-e^{-t(1+x^{2})}| \leq |t|$ uniformly for
$x\in\mathbb R$, we obtain $\|S_1(t)\| \leq 3|t|$, $\re t\geq 0$,
whenever $0\leq\beta\leq 2$. For $2<\beta\leq 4$, one has an estimate
$\|S_1(t)\| \leq M|t|e^{\re t}/(\re t)^{\beta/2-1}$  (for some $M>0$).
Hence for $\beta=4$, $S_1(t)$ is a holomorphic once integrated semigroup
(only in the extended sense if $p=\infty$), and thus Theorem \ref{Thesthol2} applies.
By using directly the estimate $\|S_1(z)\|={\cal O}(|z|/\re z)$, $\re z>0$, one obtains by integration the estimate $|z|^{\kappa-1}$ for $S_\kappa(z)$, $\kappa>2$.
\end{proof}

In conclusion, this example illustrates the full range of convergence rate for integrated
semigroups. Let us restrict ourselves to the once integrated semigroup $\{S_1(t)\}_{t\geq 0}$:
for $0\leq\beta\leq 2$, we have the optimal exponent $1$; for $2<\beta<4$
we have the exponent $2-\beta/2$, and for $\beta=4$ the operator norm $\|S_1(t)\|$ does not tend
to $0$ (the strong convergence still holds, provided $p<\infty$).
Finally for $\beta>4$ the fact that $\rho(A)=\emptyset$ makes the theory of
integrated semigroups unusable, and shows that Theorem \ref{Thesthol2} does not apply here.

\subsection{A singular semigroup with non dense range}

The following example shows that Theorem \ref{Thesthol1} can be used in more pathological cases, where the associated semigroup is not Abel summable. This is a semigroup $\{T(t)\}_{t>0}$ such that $X_0=\bigcup_{t>0}\ran T(t)$ is not dense, which implies that the subspace of strong continuity $\{x\in X: \lim_{t\rightarrow 0}\|T(t)x-x\|=0\}$ is not dense; it is inspired by \cite[Example 1.26]{Davies}.
Let $X$ be the Banach space ${\cal C}[0,1]$, and $\beta\geq 0$. We define for each $t>0$
the bounded operators $T_\beta(t)$ on $X$ by
$$
[T_\beta(t)f](x)=x^t[f(x)-f(0)(-\ln x)^\beta],\ 0<x\leq 1,\ [T_\beta(t)f](0)=0.
$$
The function $t\mapsto T_\beta(t)$ satisfies the semigroup equation and has a holomorphic
extension into the open right half-plane $\re z>0$.
If $\beta=0$, then the semigroup is degenerate. If $\beta>0$, $T(t)$ is one-to-one
but its range is not dense: in fact
$\overline{X_0}=\{f\in X: f(0)=0\}$ is the subspace of strong continuity, and the restriction of $T(t)$ to $\overline{X_0}$ is a $C_0$ semigroup.
Moreover $T(t)$ is singular when $t$ tends to $0$: one has
$$\sup_{0<x<1}|e^{z\ln x}(1-(-\ln x)^\beta)| \leq \|T(z)\| \leq 1+\sup_{0<x<1}|e^{z\ln x}(-\ln x)^\beta|,$$ 
which leads to an estimate of type $(\re z)^{-\beta}$ for the operator norm $\|T(z)\|$, $\re z>0$.

If $\beta<1$, Theorem \ref{Thesthol1} applies:
the associated $\kappa$-times integrated semigroup is bounded
for $\kappa>\beta$, admits boundary values, and satisfies the estimate (for some $M>0$)
$$
\|S_\kappa(z)\|\leq M|z|^{\kappa-\beta},\ \re z>0.
$$
The resolvent of the generator is
$$
[R(\lambda,A_\beta)f](x)=\hat{T_\beta}(\lambda)
=(\lambda-\ln x)^{-1}[f(x)-(-\ln x)^\beta f(0)].
$$
This shows that $\ran R(\lambda,A_\beta)$ is not dense in $X$ and that $\lim_{\lambda\rightarrow\infty}\lambda R(\lambda,A_\beta)f=f$ holds only for functions $f$ such that $f(0)=0$. Therefore $\{T_\beta(t)\}_{t>0}$ is not an Abel summable semigroup.

If however $\beta\geq 1$, no resolvent can be associated to $T_\beta$.
In fact the numerical function $<T_\beta(t)f,\delta_x>$ is Laplace transformable for each Dirac
measure $\delta_x\in X'$ ($x\neq 0$). Thus the resolvent, if defined, should coincide with
$$
\int_0^\infty e^{-\lambda t}<T_\beta(t)f,\delta_x>dt =
(\lambda-\ln x)^{-1}[f(x)-(-\ln x)^\beta f(0)].
$$
But this expression does not define a resolvent operator for $\beta\geq 1$:
for $\beta>1$ it is not bounded, and for $\beta=1$ the resolvent equation is not satisfied.

\subsection{Concluding remark}

In view of the special properties of hermitian integrated semigroups \cite{LiShaw},
the question arises whether there are holomorphic $\kappa$-times integrated semigroups
for $\kappa>1$ that are not obtained by integration of some holomorphic once integrated
semigroup (in the extended sense).
In other words, is it possible that the associated holomorphic semigroup satisfies
$(E^1_{\alpha,\beta})$ with $\beta>1$ and not with $\beta=1$ ?
As far as we know, the examples of such highly
singular semigroups cannot be associated to holomorphic integrated semigroups
(for example: the generator has empty resolvent set, see e.g. \cite[\S 20.5]{HP}).

\smallskip

\noindent
Department of Theoretical Physics\\
quai Ernest-Ansermet 24, \\
CH-1211 Geneva 4, Switzerland.\\
email: Vincent.Cachia@tele2.fr

\end{document}